\documentclass[11pt, a4paper]{article}  
\author{Frank Schuhmacher}

\title{Analytic decomposition of differential graded Lie algebras}

\usepackage{mathenv}
\usepackage[intlimits]{amsmath}
\usepackage[T1]{fontenc}
\usepackage{amssymb}
\usepackage{amsthm} 
\usepackage{amscd}
\usepackage[all]{xypic}
\usepackage{makeidx}

\newcounter{punkt}

\newcommand{\txbfind}[1]{\index{#1}\textbf{#1}}

\theoremstyle{definition}

\newtheorem{defi}{Definition}[section]
\newtheorem{bem}[defi]{Remark}

\newtheorem{beisp}[defi]{Example}

\newtheorem{uebung}{Exercise}

\theoremstyle{theorem}

\newtheorem{satz}[defi]{Theorem}

\newtheorem{lemma}[defi]{Lemma}
\newtheorem{kor}[defi]{Corollary}
\newtheorem{prop}[defi]{Proposition}

\newcommand{\nach}{\longrightarrow}

\newcommand{\isom}{\cong}
\newcommand{\lan}{\langle}
\newcommand{\ran}{\rangle}

\newcommand{\NN}{ \mathbb{N} }
\newcommand{\RR}{ \mathbb{R} }
\newcommand{\ZZ}{\mathbb{Z}}

\newcommand{\CC}{\mathbb{C}}

\newcommand{\KK}{\mathbb{K}}

\renewcommand{\d}{\partial}

\renewcommand{\O}{\mathcal{O}}
\newcommand{\Oh}{\mathcal{O}}

\newcommand{\m}{\mathfrak{m}}

\newcommand{\Hom}{\operatorname{Hom}}

\newcommand{\Mult}{\operatorname{Mult}}

\newcommand{\id}{\operatorname{Id}}

\newcommand{\pr}{\operatorname{pr}}

\newcommand{\Kern}{\operatorname{Kern} }

\newcommand{\ot}{\otimes}
\newcommand{\lie}{[\cdot,\cdot]}

\newcommand{\sig}{\sigma}
\newcommand{\Sh}{\operatorname{Sh}}
\newcommand{\Ot}{\operatorname{Ot}}

\newcommand{\ppp}{\cdot\ldots\cdot}
\newcommand{\ddd}{,\ldots,}
\newcommand{\kkk}{+\ldots+}
\newcommand{\odos}{\odot\ldots\odot}
\newcommand{\odots}{\odot\ldots\odot}
\newcommand{\ots}{\otimes\ldots\otimes}

\newcommand{\Wedge}{\bigwedge}
\newcommand{\pnu}{\downarrow}
\newcommand{\pno}{\uparrow}
\newcommand{\Fix}{\operatorname{Fix}}

\newcommand{\Der}{\operatorname{Der}}
\newcommand{\Coder}{\operatorname{Coder}}

\newcommand{\DGmanf}{DG-\textsc{Manf}}

\newcommand{\norm}{||\cdot||}
\newcommand{\eps}{\epsilon}
\newcommand{\lam}{\lambda}
\newcommand{\CODER}{\mathcal{C}\! \it{oder}}

\newcommand{\Vect}{\operatorname{Vect}}

\newcommand{\hQ}{\hat{Q}}

\begin{document}
\maketitle

\begin{abstract}
We prove explit formulas for the decomposition of 
a differential graded Lie algebra
into a minimal and a linear $L_\infty$-algebra.
We define  a
category of metric $L_\infty$-algebras,
called Palamodov $L_\infty$ algebras, where the structure
maps satisfy a certain convergence condition
and deduce a decomposition theorem 
for differential graded Lie algebras
in this category. This theorem serves for instance
to prove the convergence of the Kuranishi map
assigned to a differential graded Lie algeba.
\end{abstract}

\tableofcontents


\section*{Introduction}

As observed by Kontsevich \cite{KontDQ}
and others, any $L_\infty$-algebra over a field
of characteristic zero can be
decomposed into a direct sum of a minimal 
and a linear contractible $L_\infty$-algebra.
This means that
there exists an $L_\infty$-structure on the
cohomology $H(L)$ of the complex $(L,d)$,
where $d$ is the linear component of the 
$L_\infty$-structure of $L$, such that $H(L)$ is
a direct summand of $L$, not only as a
graded module but also as an $L_\infty$-algebra.
This has two important consequences: first,
any quasi-isomorphism between $L_\infty$-algebras
over a field of characteristic zero can
be inverted (we may thus speak of $L_\infty$\textit{-equivalence})
and secondly, the
deformation functors
(local Artin rings)$\to$(sets) 
assigned to quasi-isomorphic DG Lie algebras
(see \cite{KontDQ} or \cite{Mane})
are isomorphic. Thus, if the deformations of an object are
\textit{governed} by a DG Lie algebra $L$ (i.e. deformations of
the object correspond to Maurer-Cartan elements of $L$
and equivalent deformations are in the same orbit
of the action of a Lie group related to $L^0$),
then the 
$L_\infty$-equivalence class of $L$ contains the
formal deformation theory of the object.
For example, the base of a formally versal
deformation, if it exists, is defined by the
\textit{Kuranishi-map} $H(L)^1\nach H(L)^2$,
whose Taylor coefficients are the
restrictions of the $L_\infty$-structure maps 
of $H(L)$.
This motivates the question
of an explicit description of the $L_\infty$-structure
on $H(L)$ and of an $L_\infty$-map $H(L)\to L$
in the case where $L$ is just a DG Lie algebra
(recursion formulas were already indicated by 
Kontsevich \cite{KontL} and Huebschmann/Stasheff \cite{Hust}).
We give an explicit description for the structure on $H(L)$
in function of a splitting
(i.e. a degree $-1$ map $\eta:L\to L$ with $d\eta d=d$)
of the complex $(L,d)$, and an explicit description
of an $L_\infty$-isomorphism $H(L)\oplus (F,d)\to L$, 
where $F$ is the complement
of $H(L)$ in $L$. The improvement with respect to
\cite{Schuh} is that now, we also have an explicit
description of the map $(F,d)\to L$. In \cite{Schuh},
the existence of such a map was obtained by
obstruction theory. 
Since an inverse isomorphism can 
recursively be constructed by the inverse function theorem
for $L_\infty$-algebras,
this allows a complete algorithmic decomposition of
a split DG Lie algebra.\\ 

There is the deeper question of
the convergence of the Kuranishi map.
The present paper contains a quite general answer
to this question: To formulate the
main result, we need some analytic language:
First of all, the algebraic language of
$L_\infty$-algebras or equivalently, of formal DG manifolds,
is not sufficient for
convergence questions. It would be straight
forward, but for many deformation problems
still not sufficient, to define
Banach $L_\infty$-algebras, i.e. $L_\infty$-structures
$(\mu_n)_{n\in\NN}$
on underlying graded Banach spaces $L$, such that
the formal sums $\sum_n\frac{1}{n!}\tilde\mu_n$,
with $\tilde\mu_n(x)=\mu_n(x\ddd x)$ 
converge on some zero neighbourhood.
We introduce, more generally,
\textit{Palamodov $L_\infty$-algebras,} where
the underlying graded vector spaces do not
even have a norm but a metric defined by a countable
family of semi-norms.
We need a subtable notion of convergence, for this context.
For \textit{Palamodov DG Lie algebras} admitting
a splitting, we prove the convergence
of the structure 
$\mu$ on $H(L)$
and of the map $H(L)\oplus F\to L$ and
since the inverse function theorem is valid
for Palamodov $L_\infty$-algebras, this implies 
the decomposition 
$$(L,d,\lie)\isom (H,\mu)\oplus (F,d)$$
in the category of Palamodov $L_\infty$-algebras.
Recall that this decomposition specialises to the
Banach context.
\pagebreak

\section{Formal DG manifolds}\label{general} 

In this section, we give a detailed review
of the correspondence between $L_\infty$-algebras,
formal DG-manifolds and free DG coalgebras.
Experts may skip this section. 

\subsection{Graded symmetric and exterior algebras}

Let $\KK$ be a commutative ring of characteristic zero.
For a graded $\KK$-module $M$, 
the graded symmetric algebra $SM$ is
defined as the tensor algebra $TM=\oplus_{n\geq 0}M^{\ot n}$
modulo the relations $m_1\ot m_2-(-1)^{m_1 m_2}m_2\ot m_1=0$.
We denote the graded symmetrical product by $\odot$.
The algebra $TM$ (resp. $SM$) is bigraded. 
The graduation on $TM$ (resp. $SM$) defined by 
$g(m_1\ot\ldots \ot m_n)=g(m_1)+\ldots +g(m_n)$ 
(resp. $g(m_1\odot\ldots \odot m_n)=g(m_1)+\ldots+ g(m_n)$), 
where $g$ is the graduation of $M$, will be called 
\txbfind{linear graduation}.
The one defined by $g(m_1\ot\ldots \ot m_n)=n$ 
(resp. $g(m_1\odot\ldots \odot m_n)=n$) 
will be called \txbfind{polynomial graduation}.
On $S_+M:=\oplus_{n\geq 1}M^{\odot n}$,\index{$S_+M$}
there is a natural $\KK$-linear comultiplication
$\Delta^+:S_+M\to S_+M\ot S_+M$, given by
$$m_1\odot\ldots \odot m_n\mapsto\sum_{j=1}^{n-1}\sum_{\sig\in\Sh(j,n)} 
\epsilon(\sigma,m_1,\ldots ,m_n)m_{\sig(1)}\odos m_{\sig(j)}\ot
m_{\sig(j+1)}\odos m_{\sig(n)},$$
where the sign $\epsilon(\sig):=\epsilon(\sigma,m_1,\ldots ,m_n)$ 
is defined in such a way that
$m_{\sig(1)}\odos m_{\sig(n)}=\epsilon(\sig)m_1\odos m_n$.
Note that $\Kern\Delta^+=M$.
On $SM$, there is a $\KK$-linear comultiplication $\Delta$, 
defined by $\Delta(1):=1\ot 1$ and 
$\Delta(m):=m\ot 1+\Delta'(w)+1\ot m$,
for $m\in S_+M$. Note that $\Delta$ is injective.

For a graded module $L$, the graded exterior algebra $\Wedge^+ L$ 
without unit is\index{graded exterior algebra}\index{$\Wedge^+$}
defined as the tensor algebra $T_+L=\oplus_{n\geq 1}L^{\ot n}$
modulo the relations $a_1\ot a_2+(-1)^{a_1 a_2}a_2\ot a_1=0$.
We denote the graded exterior product by $\wedge$.
By $L[1]$, we denote the graded module with $L[1]^i=L^{i+1}$
and by $\pnu$ the canonical map $L\to L[1]$ of degree $-1$.
Set $\pno:=\pnu^{-1}$. \index{$\pnu$, $\pnu$}
Remark that for each $n\geq 1$, there is an isomorphism
\begin{align*}
\pnu^n:{\Wedge}^n L&\nach L[1]^{\odot n}\\
a_1\wedge\ldots \wedge a_n & \mapsto (-1)^{(n-1)\cdot a_1\kkk 1\cdot\ldots
  a_{n-1}} \pnu a_1\odos\pnu a_n.
\end{align*}
Its inverse map is given by $(-1)^{\frac{n(n-1)}{2}}\pno^n$.
In this formula, we deduce the sign from
the Koszul convention. More generally, 
for homogeneous graded morphisms $f,g$
of graded modules, we set $(f\ot g)(a\ot b):=(-1)^{ga}f(a)\ot g(b)$. 
In the exponent, $a$ always means the degree
of an homogeneous element (or morphism) $a$ and $ab$ means
the product of degrees and not the degree of the product.\\

For $\sig\in\Sigma_n$ and $a_1\ddd a_n\in L$, we define the
sign $\chi(\sig):=\chi(\sig,a_n\ddd a_n)$ in such a way that
$$a_{\sig(1)}\wedge\ldots \wedge a_{\sig(n)}=\chi(\sig)a_1\wedge
\ldots \wedge a_n.$$
The 
correlation between $\chi$ and $\epsilon$ is given by
$$\chi(\sig,a_1\ddd a_n)=
(-1)^{(n-1)(a_1+a_{\sig(1)})\kkk 1\cdot(a_{n-1}+a_{\sig(n-1)})}
\epsilon(\sig,\pnu a_1\ddd\pnu a_n),$$
for $a_1,\ldots ,a_n\in L$.
For a graded module $V$,
we define two different actions of the symmetric group
$\Sigma_n$ on $V^{\ot n}$:
The first one is given by
$$
(\sig,v_1\ot\ldots \ot v_n)\mapsto\epsilon(\sig,v_1,\ldots ,v_n)
v_{\sig(1)}\ot\ldots \ot v_{\sig(n)}.$$
The application of a $\sig$
commutes with the canonical projection
$V^{\ot n}\to V^{\odot n}$.
The second one is given by
$$(\sig,v_1\ot\ldots \ot v_n)\mapsto\chi(\sig,v_1,\ldots ,v_n)
v_{\sig(1)}\ot\ldots \ot v_{\sig(n)}.$$
Here, the application of a $\sig$
commutes with the canonical projection
$V^{\ot n}\to \wedge^nV$.
When we work with symmetric powers, we use the first
action; when we work with exterior powers, we use the
second one. Since the context shall always be clear, 
we don't distinguish both actions by different notation.
We will use the (anti-)symmetrization maps:\index{$\alpha_n$}
$$\alpha_n:=\sum_{\sig\in\Sigma_n}\sig.:V^{\ot n}\nach V^{\ot n}.$$
When $\sig.$ denotes the first action, $\alpha_n$
can be considered as a map $V^{\odot n}\to V^{\ot n}$; when $\sig.$
denotes the second action, $\alpha_n$
can be considered as a map $\wedge^n V\to V^{\ot n}$.
Furthermore, for both cases, we will use the maps
$$\alpha_{k,n}:=\sum_{\sig\in\Sh(k,n)}\sig.:
V^{\ot n}\nach V^{\ot n}.$$
For the natural projection $\pi:V^{\ot n}\to V^{\odot n}$
(resp. $V^{\ot n}\to \Wedge^nV$), we have
$$\pi\circ\alpha=n!\id.$$

\subsection{Free differential graded coalgebras}

Let $(C_1,\Delta_1)$ and $(C_2,\Delta_2)$ be coalgebras.
Remember that a module homomorphism 
$F:C_1\to C_2$ is a coalgebra morphism,
if and only if the diagram
\index{coalgebra morphism}
\begin{equation}\label{qudiag}
\xymatrix{
C_1\ar[r]^F\ar[d]^{\Delta_1} & C_2\ar[d]^{\Delta_2}\\
C_1\ot C_1\ar[r]^{F\ot F} & C_2\ot C_2
}\end{equation}
commutes. Each coalgebra morphism
$F:(SM,\Delta)\to (SM',\Delta)$ satisfies $F(1)=1$. 
The restriction $F\mapsto F|_{S_+M}$
is a one-to-one correspondence
between coalgebra morphisms
$(SM,\Delta)\to (SM',\Delta)$ and coalgebra morphisms
$F\!:(S_+M,\Delta^+)\! \to (S_+M',\Delta^+)$.
We fix the following notions:

\begin{align*}
 \hat{F}_n:=&  F|_{M^{\odot n}}: M^{\odot n}\nach  SM'\\
F_{k,l}:=& \pr_{{M'}^{\odot l}}\circ\hat{F}_k: M^{\odot k}
\nach {M'}^{\odot l}\\
F_n:=&F_{n,1}: M^{\odot n}\nach M'
\end{align*}

Sometimes, we shall consider the maps $F_n$ as graded symmetric
maps $M^{\ot n}\to M'$ instead of maps $M^{\odot n}\to M'$. 
For each multi-index $I=(i_1,\ldots ,i_k)\in\NN^k$, we set
$I!:=i_1!\ppp i_k!$ and $|I|:=i_1\kkk i_k$ and
$$F_I:=\frac{1}{I!k!}(F_{i_1}\odot\ldots \odot F_{i_k})\circ\alpha_n.$$
Here, by $F_{i_1}\odot\ldots \odot F_{i_k}$, we mean the composition
of $F_{i_1}\ot\ldots \ot F_{i_k}$ and the natural projection
${M'}^{\ot k}\to {M'}^{\odot k}$.
A morphism $F:SM\to SM$ of graded 
coalgebras is called \txbfind{strict},
if $F_n=0$ for each $n\geq 2$.\\ 

\index{coderivation}
For a coalgebra $(C,\Delta)$,
remember that a module homomorphism $Q:C\to C$ is a coderivation,
if and only if the diagram
\begin{equation}\label{qdelta}
\xymatrix{
C\ar[rr]^{Q}\ar[d]^\Delta && C\ar[d]^{\Delta}\\
C\ot C\ar[rr]^{Q\ot 1+1\ot Q}& & C\ot C
}\end{equation}
commutes.
We fix the following notions:

\begin{align*}
\hQ_n:=&Q|_{M^{\odot n}}: M^{\odot n}\nach SM\\
Q_{k,l}:=&\pr_{M^{\odot l}}\circ \hQ_k: M^{\odot k}\nach M^{\odot l}\\
Q_n:=&Q_{n,1}: M^{\odot n}\nach M
\end{align*}

Observe that there is a one-to-one correspondence between
degree $i$ coderivations on $(S_+M,\Delta^+)$
and degree $i$ coderivations of $Q$ degree $+1$ 
on $(SM,\Delta)$ with $Q(1)=0$. 
By the following proposition, there is a one-to-one correspondence
between degree $i$ coderivations $Q:SM\to SM$ and 
sequences of linear maps $Q_n:S_nM\to M$ of degree $i$, and
a one-to-one correspondence
between coalgebra maps $ F:SM\to SM'$
and sequences of linear maps 
$F_n:S_nM\to M'$, $n\geq 1$.

\begin{prop}
Let $Q$ (resp. $Q'$) be a degree $i$ coderivation 
on $SM$ (resp. $SM'$) and 
$F:=SM\to SM'$ a morphism of graded coalgebras. 
Then, for $n\geq 1$ and
$1\leq l\leq n+1$, $1\leq k\leq n$,  we have
\begin{equation*}
Q_{n,l}=
(Q_{n-l+1}\ot 1\ot\ldots \ot 1)\circ\alpha_{n-l+1,n}
\end{equation*}
and
\begin{equation*}
F_{n,k}=\sum_{i_1\kkk i_k=n}F_I.
\end{equation*}
The map $F$ respects the coderivations $Q$ and $Q'$ if and only if
$F(Q(1))=Q'(1)$ and for each $n\geq 1$,
we have
\begin{equation}\label{coalgmor}
\sum_{k=1}^n\sum_{I\in\NN^k\atop |I|=n}
Q'_k\circ F_I=
\sum_{k+l=n+1}
F_l\circ(Q_{k}\ot 1\ots 1)\circ\alpha_{k,n}.
\end{equation}
On the right hand-side, the sum is over all $l\geq1$ and $k\geq 0$.
The term
$(Q_0\ot 1\ots 1)(m_1\ots m_n)$ must be interpreted as
$Q_0(1)\ot m_1\ots m_n$.
\end{prop}

Set $\Coder^i(SM)$
to be the module of degree $i$ coderivations.
The direct sum $\Coder(SM)=\oplus_{i\in\ZZ}\Coder^i(SM)$
is a graded Lie algebra with graded commutator
$$[Q,q]=Q\circ q-(-1)^{Q\cdot q}q\circ Q$$
as Lie bracket. Explicitly,
\begin{equation}\label{explizit}
[Q,q]_n=
\sum_{k+l=n+1} Q_l\circ (q_k\ot 1^{\ot n-k})\circ\alpha_{k,n}-
(-1)^{Q\cdot q} q_l\circ (Q_k\ot 1^{\ot n-k})\circ\alpha_{k,n}.
\end{equation}
Observe that a degree $1$ coderivation
$Q$ on $SM$ is a codifferential (i.e.
$Q^2=0$) if and only if, for each $n\geq 0$, we have 
$$\sum_{k+l=n+1}Q_l(Q_k\ot 1^{\ot n-k})\circ\alpha_{k,n}=0.$$
By definition, a free DG coalgebra is
a pair $(SM,Q^M)$, where $M$ is a graded $\KK$-module and
$Q^M$ is a degree one
codifferential on $SM$. A morphism of free
DG coalgebras is a coalgebra map of degree zero
respecting the codifferentials.

\subsection{Formal vector fields as coderivations}

Before we introduce the equivalent notions of formal
DG manifolds and $L_\infty$-algebras, we make some
remarks on analytic space germs and on analytic vector fields:
An analytic function $f:\CC^N\nach\CC^m$ defined on a
neighbourhood of $0$ can weather
be represented by a power series
$$f(x)=\sum_\nu a_\nu x^\nu$$
with coefficients $a_\nu$ in $\CC^m$
or as a sum
$$f(x)=\sum_{n=0}^\infty\frac{1}{n!}\tilde{f}_n(x),$$
where $\tilde f_n(x)=f_n(x\ddd x)$, for a
uniquely defined symmetric $n$-linear map 
$f_n:\CC^N\ots\CC^N\nach\CC^m$.
The map $\tilde f_n$ is the $n$-th Taylor coefficient of $f$. 
The relation between the coefficients $a_\nu$ and $f_n$ is
given by the formulas
$$f_n=\sum_{|\nu|=n}a_\nu (e^\ast)^{\odot\nu}$$
and
$$a_\nu=\frac{1}{n!}f_{|\nu|}(e^{\odot\nu}).$$
For the verification of this formulas, remark that
$$(a_\nu (e^\ast)^{\odot\nu})^\sim (x)=
\frac{|\nu|!}{\nu!}a_\nu(e^\ast)^{\odot\nu}(x^{\odot\nu})=
|\nu|!a_\nu x^\nu.$$

Let $R=\CC\lan e_j\ran_{j\in J}$ be a free g-finitely generated 
algebra over the convergent power series algebra 
$R^0=\CC\{e_j\}_{j\in J^0}$. Let $\hat R=\CC [[x_j]]$ be the $\m$-adic
completion
of $R$,  where $\m$ is the ideal generated by the $x_j;\;j\in J$.
Consider the graded $\CC$-vector space $M=\coprod_{j\in J}\CC e_j$
with generators of degree $-g(e_j)=g(x_j)$.
Consider $SM$ as graded coalgebra.
The multiplication $\odot$ in the dual algebra $(SM)^\ast$ of 
$(SM,\Delta)$ is given by
$$(a\odot b)(m)=(a\ot b)(\Delta(m))
=\sum_{\sig\in\Sh(k,k+l)}\eps(\sig,m)(-1)^{m'_\sig\cdot b}a(m'_\sig)\cdot
b(m''_\sig),$$
for $a\in(M^{\odot k})^\ast$, $b\in (M^{\odot l})^\ast$
and $m\in M^{\odot k+l}$.
We have an isomorphism $\hat R\to (SM)^\ast$
of graded algebras sending $x_j$ to
the dual generator $e_j^\ast$ of $M^\ast$.
The inverse map is given by
$$f\mapsto\sum_{n=0}^\infty\sum_{|\nu|=n}\frac{1}{\nu!}
f_n(e^{\odot\nu})\cdot x^\nu,$$
where we sum only over the ordered multi-indices with indices
of odd degree appearing at most once.
For an ordered multi-index $\nu=(\nu_1\ddd\nu_N)$,
we have
$$(e^\ast)^{\odot\nu}:=(e_1^\ast)^{\odot\nu_1}\odots
(e_N^\ast)^{\odot\nu_N}=\nu!(e^{\odot\nu})^\ast.$$

Let $\m$ be the canonical maximal
ideal of $\hat R$ and set 
$\Der(\hat R)=\oplus_{i\in\ZZ}\Der^i(\hat R)$
to be the graded vector space of all $\CC$-linear,
$\m$-adic continuous derivations of $\hat R$.
The graded commutator
defines a graded Lie structure on $\Der(\hat R)$.
Observe that the dualization map $Q\mapsto Q^\ast$
is an isomorphism of graded Lie algebras
$$\Coder(SM)\nach\Der(S(M)^\ast).$$
The inverse map sends a derivation
$\delta$ to the sum
$\sum_{j\in J}e_j\cdot\delta(e_j^\ast)$.\\

Consider an open subset $M$ of
$\CC^m$ with $0\in M$. An analytic vectorfield $q$ on $M$ is 
an analytic map $M\to\CC^m$. Locally at $0$, it is given by a power 
series
$$\tilde q:x\mapsto\sum_{n=0}^\infty \frac{1}{n!}q_n(x\ddd x),$$
where $q_n:\CC^m\times\ldots\times\CC^m\to\CC^m$ is a unique
symmetric n-linear map, i.e. $\tilde q_n:x\mapsto q_n(x\ddd x)$
is a homogeneous polynomial function of degree $n$.
Recall that as a Lie algebra, we can identify the set of
analytic vectorfields on $M$ with the Lie algebra of $\CC$-linear
derivations of the section algebra $\Gamma(\O_M,M)$, where
the Lie bracket is the commutator of derivations. The following
exercise gives the algebraic description of this Lie bracket:

\begin{uebung}
For vectorfields $Q,q$ on $M$, we have 
$[Q,q](x)=\sum_{n=0}^\infty \frac{1}{n!}[Q,q]^\sim_n (x)$, where the
symmetric n-linear map $[Q,q]_n$ is given by the following formula:
\begin{align}\label{klammer}
\begin{array}{c}
[Q,q]_n(x_1\ddd x_n)=\hfill\\
\quad\quad\sum_{k+l=n+1}\sum_{\sig\in\Sh(k,n)}
Q_l(q_k(x_{\sig(1)}\ddd x_{\sig(k)}),x_{\sig(k+1)}\ddd x_{\sig(n)})-\\
\hfill q_l(Q_k(x_{\sig(1)}\ddd x_{\sig(k)}),x_{\sig(k+1)}\ddd x_{\sig(n)}).
\end{array}
\end{align}
\end{uebung}
If we replace $M$ by an arbitrary vector space,
we can define the Lie algebra of \textbf{formal
vectorfields} on $M$, i.e. of formal sums 
$q=\sum\frac{1}{n!}\tilde q_n$, where
each $q_n:M^{\ot n}\to M$ is a symmetric n-linear map and
the bracket $[Q,q]$ of formal vectorfields $Q,q$ on $M$
is the formal sum $\sum [Q,q]_n$, where $[Q,q]_n$ is defined
by the equation~(\ref{klammer}).
To define formal DG manifolds, we generalise the definition 
of formal vectorfields to
the graded context, where the word \textit{symmetry} is replaced 
by the word \textit{graded symmetry}.

\subsection{Formal DG manifolds and  $L_\infty$-algebras}

Let $M$ be a graded $\KK$-module. A \textbf{formal
vectorfield of degree $i\in\ZZ$} on $M$ is
a formal sum 
$q=\sum\frac{1}{n!}\tilde q_n$, where
each $q_n:M^{\ot n}\to M$ is a graded symmetric n-linear map
of linear degree $i$. Set $\Vect^i(M)$ to be the $\KK$-module
of formal vectorfields of degree $i$ on $M$.
The direct sum $\Vect(M):=\sum_{i\in\ZZ}\Vect^i(M)$ has
a graded Lie algebra structure with Lie
bracket defined by equation (\ref{explizit}).
By definition, the graded Lie algebras $\Vect(M)$
and $\Coder(SM)$ are canonically isomorphic.
A \textbf{formal DG manifold} is a pair $(M,Q^M)$,
where $M$ is a graded $\KK$-module and $Q^M$
is a formal degree one vectorfield on $M$ with
$[Q^M,Q^M]=0$. In other words, a formal DG manifold
is a pair $(M,Q^M)$ such that the pair
$(SM,Q^M)$ is a free DG coalgebra. A \textbf{morphism}
$f:M\nach M'$ of formal DG manifolds is a formal
sum $f=\sum_{n=0}^\infty\tilde f_n$, where
$f_n:M^{\ot n}nach M$ is a graded symmetric map such
that the family $(f_n)_{n\geq 1}$ defines a morphism
$SM\nach SM'$ of DG coalgebras.
Thus, the category $\DGmanf$ of formal
DG manifolds is canonically isomorphic to
the category of free DG coalgebras.
Geometrically, a morphism of formal DG manifolds is
a morphism of space germs, equivariant with respect
to the given vectorfields.\\

Remember that a module $L$ with a sequence of maps 
$\mu_n:\Wedge^nL\nach L$ of degree $2-n$, for $n\geq 0$, is called an
\txbfind{$L_\infty$-algebra} if the 
pair $(M,Q^M)$ with
$M:=L[1]$, and 
$$Q_n:=(-1)^{n(n-1)/2}\pnu\circ\mu_n\circ\pno^n:M^{\odot n}\nach M$$ 
is a formal DG manifold.
This condition is equivalent to the equations
\begin{align}
\sum_{k+l=n+1}\sum_{\sig\in\Sh(k,n)}(-1)^{k(l-1)}\chi(\sig)
\mu_l(\mu_{k}(a_{\sig(1)}\ddd a_{\sig(k)}), a_{\sig(k+1)}\ddd 
a_{\sig(n)})=0
\end{align}for each $n\geq 0$ and $a_1\ddd a_n\in L$.
In the literature, $\mu_0$ is mostly assumed to be trivial.
If this is the case, $(L,\mu_1)$ is a DG module. 
An $L_\infty$-algebra $(L,\mu_n)_{n\geq 1}$ is called \txbfind{minimal}, if
$\mu_1=0$. It is called \txbfind{linear}, if
$\mu_i=0$ for $i\geq 2$.\\

Let $(L,\mu_\ast)$ and $(L',\mu'_\ast)$ be $L_\infty$-algebras.
Set $M:=L[1]$, $M':=L'[1]$.
A sequence of maps
$f_n:\Wedge^nL\nach L';\;n\geq 1$ of degree $1-n$ is called
\txbfind{$L_\infty$-morphism}, if the maps $F_n:=W^{\odot n}\nach W'$ 
induced by $f_n$ (explicitly: 
$F_n=(-1)^{n(n-1)/2}\pnu\circ f_n\circ\pno^n$) define a morphism
of formal DG manifolds. 
Rewrite
condition (\ref{coalgmor}) in terms of $f_n$ and $\mu_n$:

\begin{equation}\label{beding}
\begin{array}{l}
\mu'_1\circ f_n-\sum_{i+j=n}
\frac{(-1)^i}{2}\mu'_2(f_i,f_j)\circ\alpha_{i,j}\\[4mm]
+\quad\sum_{k=3}^n\sum_{I\in\NN^k\atop |I|=n}
(-1)^{k(k-1)/2+i_1(k-1)+\ldots +i_{k-1}\cdot 1}\mu'_k\circ f_I\\[3mm]
\quad=\quad
\sum_{k+l=n+1}(-1)^{k(l-1)}f_l\circ(\mu_k\ot 1\ot\ldots \ot 1)\circ\alpha_{k,n}
\end{array}
\end{equation}

\noindent
For the case
where $L'$ is a differential graded Lie-algebra, i.e.
$\mu'_k=0$ for $k=0$ and $k\geq 3$, set
$d:=\mu'_1$ and $[\;,\;]:=\mu'_2$. We get the following
conditions for the maps $f_n$ to define an $L_\infty$-morphism
(see Definition 5.2 of \cite{Lama}):

$$\begin{array}{l}
df_n(a_1\ddd a_n)-\\[3mm]
\sum_{i+j=n}\sum_\sigma\chi(\sig)
(-1)^{i+(j-1)(a_{\sig(1)}\kkk a_{\sig(i)})}[f_i(a_{\sig(1)}\ddd
a_{\sig(i)}),f_j(a_{\sig(i+1)}\ddd a_{\sig(n)})]\\[4mm]
\quad=\sum_{k+l=n+1}\sum_{\sig\in\Sh(k,n)}
(-1)^{k(l-1)}\chi(\sig)f_l(\mu_k(a_{\sig(1)}\ddd
a_{\sig(k)}),a_{\sig(k+1)}\ddd a_{\sig(n)}),
\end{array}$$
where $a_1,\ldots ,a_n\in L$ and
the second sum goes over all $\sig$ in $\Sh(i,n)$
such that $\sig(1)<\sig(i+1)$.
A morphism $f:(L,\mu_n)_{n\geq 1}\nach (L',\mu'_n)_{n\geq 1}$ of 
$L_\infty$-algebras is called an
\txbfind{$L_\infty$-quasi-isomorphism},
if $f_1$ is a quasi-isomorphism of differential graded modules.

\begin{prop}\label{isom}
Let $f:M\nach M'$ be a morphism of formal supermanifolds.
Suppose, there is a module homomorphism $g':M'\nach M$ such that
$g'\circ f_1=\id_M$. Then, there is a morphism $g:M'\nach M$
of formal supermanifolds such that $g_1=g'$ and $gf=\id_M$.
If $f_1$ is an isomorphism with inverse $g$ and if
$f$ is $Q$-equivariant and if
$g'$ respects $Q_1^{M'}$ and $Q_1^M$, then $g$ can be chosen 
$Q$-equivariant as well.
\end{prop}
\begin{proof}
One can check directly that the 
sequence of maps defined by
$$g_n:=-\sum_{k=2}^n\sum_{I\in\NN^k\atop |I|=n}g_1\circ f_k\circ (g_{I}),$$
for $n\geq 2$, define a morphism of formal supermanifolds
with the desired property.
\end{proof}

\section{Trees}\label{trees}

Trees were first used by
Kontsevich/ Soibelman \cite{KontS} to describe the
$A_\infty$-structure that a DG module, homotopy equivalent to
a differential graded algebra inherits. 
We have a similar objective, but for $L_\infty$-algebras
instead of $A_\infty$-algebras.
Trees in our definition are always binary trees. 
We give a definition of binary trees 
and assign several invariants to them, which
are important in order to get good signs later on.

\subsection{Definitions}

A \textbf{(binary) tree} with $n$ leaves 
consists of a pair $\phi=(\phi,V)$ where
$V=\{K_0\ddd K_{n-2} \}$ denotes a set
of \textbf{ramifications} 
and $\phi:\{K_1,\ldots,K_{n-1}\}\nach\{K_0,\ldots,K_{n-1}\}$
is a map, such
that for each $i=0\ddd n-2$, 
\begin{enumerate}
\item[(i)]
the inverse image $\phi^{-1}(K_i)$ contains at most 2 elements;
\item[(ii)]
there is an $n\geq 0$ such that $\phi^n(K_i)=K_0$.
\end{enumerate}
$K_0$ is called \textbf{root} of $\phi$. 
There is a tree with one leaf and no ramification
which will always be denoted by $\tau$.
\index{$\tau$}
An \textbf{orientation} of a tree $(\phi,V)$ is a family
$\pi=(\pi_K)_{K\in V}$ of inclusions $\pi_K:\phi^{-1}(K)\nach\{1,2\}$.
The triple $\phi=(\phi,V,\pi)$ is called an \textbf{oriented tree.}
For each oriented tree $(\phi,V,\pi)$, there is a natural 
\textbf{ordering}
on the set $V$: For $K\in V\setminus K_0$, 
suppose that $\phi^m(K)=K_0$.
We set
$$v(K):=\frac{\pi_{\phi(K)}(K)}{3^m}+
\frac{\pi_{\phi^2(K)}(\phi(K))}{3^{m-1}}\kkk
\frac{\pi_{\phi^m(K)}(\phi^{m-1}(K))}{3}.$$
Set $v(K_0):=0$. Then $v:V\nach\RR$ is injective, hence it induces
an ordering on $V$.
When we write down the value $v(K)$ of a ramification $K$ in its
3-ary decomposition, we just get an algorithm, how to get 
from the root $K_0$ 
to $K$. For example $0.1121$ means 
``go (in the driving direction) right-right-left-right''. 
When $(\phi,V,\pi)$ is an oriented tree with $n$ leaves, we can extend
the map $\phi$ to a map 
$\tilde{\phi}:V\setminus{K_0}\cup\{1,\ldots ,n\}\nach V$ such
that
\begin{enumerate}
\item[(i)]
For $1\leq i< j\leq n$, we have $\tilde{\phi}(i)\leq\tilde{\phi}(j)$.
\item[(ii)]
For each $K\in V$, $\tilde{\phi}^{-1}(K)$ has exactly 2 elements.
\end{enumerate}
 
The numbers $1,\ldots ,n$ stand for the leaves of $\phi$.
Furthermore, we can extend the map $v:K\nach [0,1)$ on
$\tilde{V}:=V\cup\{1,\ldots ,n\}$ in such a way that the 3-ary
decomposition of $v(i)$ describes the way from the root to the
i-th leaf of $\phi$, for $i=1,\ldots ,n$. Then we have
$v(i)<v(j)$, for $1\leq i<j\leq n$. In consequence, we have
an ordering on $\tilde{V}$.\\

Two trees $(\phi,V)$ and $(\phi',V')$ are called \textbf{equivalent},
if there is a bijection $f:V\nach V'$ of the ramification sets
such that $f\circ\phi=\phi'\circ f$.
Two oriented trees $(\phi,V,\pi)$ and $(\phi',V',\pi')$ are called
\textbf{oriented equivalent},
if there is a bijection $f:V\nach V'$ of the ramification sets
such that $f\circ\phi=\phi'\circ f$ and $\pi'\circ f=\pi$.  
When we draw oriented trees, we shall put elements $K'$ of 
$\phi^{-1}(K)$ down left of $K$ if $\pi_K(K')=1$ and down right 
of $K$ if $\pi_K(K')=2$.

\begin{beisp}\label{extree}
The following trees with three leaves are equivalent
but not oriented equivalent: 

{\footnotesize
\setlength{\unitlength}{0.4cm}
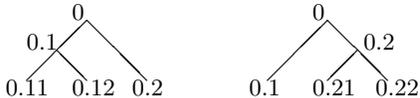
\begin{picture}(14,4)
\put(1,1){\line(1,1){2}}
\put(2,2){\line(1,-1){1}}
\put(3,3){\line(1,-1){2}}
\put(9,1){\line(1,1){2}}
\put(11,3){\line(1,-1){2}}
\put(11,1){\line(1,1){1}}
\put(0.3,0.5){$0.11$}
\put(2.5,0.5){$0.12$}
\put(1,2){$0.1$}
\put(2.5,3){$0$}
\put(8.4,0.5){$0.1$}
\put(4.5,0.5){$0.2$} 
\put(10.5,3){$0$}
\put(10.5,0.5){$0.21$}
\put(12.2,2){$0.2$}
\put(12.6,0.5){$0.22$}
\end{picture}
}

For each ramification and each leaf, we have indicated its value.
\end{beisp}
\index{$\Ot(n)$}
Set $\Ot(n)$ to be the set of equivalence classes of 
oriented trees with
$n$ leaves. Observe that the cardinality $\#\Ot(n)$ of $\Ot(n)$
is given by the
$(n-1)$-st Catalan number $c_{n-1}=\frac1n {2(n-1) \choose n-1}$
which equals $\frac{(2(n-1))!}{n!(n-1)!}$.
Using Stirling's
$$\sqrt{2\pi n}\cdot (\frac ne)^n\leq n!\leq2\cdot
\sqrt{2\pi n}\cdot(\frac ne)^n,$$
we get the following:
\begin{lemma}\label{otnum}
$\#\Ot(n)<4^{n-1}$.
\end{lemma}
The set $\Ot(2)$ contains just one element. We denote it by
$\beta$.
The set $\Ot(4)$ contains just the following elements:

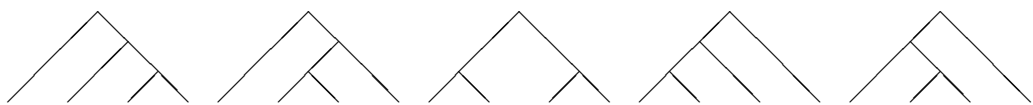
\setlength{\unitlength}{0.4cm}
\begin{picture}(33,5)
\put(1,1){\line(1,1){3}}
\put(3,1){\line(1,1){2}}
\put(5,1){\line(1,1){1}}
\put(4,4){\line(1,-1){3}}

\put(8,1){\line(1,1){3}}
\put(10,1){\line(1,1){2}}
\put(11,2){\line(1,-1){1}}
\put(11,4){\line(1,-1){3}}

\put(15,1){\line(1,1){3}}
\put(19,1){\line(1,1){1}}
\put(16,2){\line(1,-1){1}}
\put(18,4){\line(1,-1){3}}

\put(22,1){\line(1,1){3}}
\put(23,2){\line(1,-1){1}}
\put(24,3){\line(1,-1){2}}
\put(25,4){\line(1,-1){3}}

\put(29,1){\line(1,1){3}}
\put(31,1){\line(1,1){1}}
\put(31,3){\line(1,-1){2}}
\put(32,4){\line(1,-1){3}}
\end{picture}

For a tree $(\phi,V)$ and $K\in V$, there is a tree $\phi|_K$
with root $K$ and ramifications
$\{ K'\in V:\;\phi^n(K')=K \text{ for an } n\geq 0\}$.
We have to introduce several invariants:
For a tree $\phi$ with $n>1$ leaves and $1\leq i\leq n$, set $w_\phi(i)$
to be the difference of the number $s_\phi(i)$ of ramifications
of $\phi$ which are smaller than $i$ and $i-1$. ($i-1$ is the
number of leaves of $\phi$, smaller than $i$.) 
For $K\in V$, set $w_\phi(K):=w_{\phi-\phi|_K}(K)$, where on the right
hand-side, $K$ is considered as leaf of $\phi-\phi|_K$.
\index{$s_\phi$}\index{$w_\phi$}\index{$e(\phi)$}

\begin{bem}
For $K\in V$, the integer $w_\phi(K)$ is just the number of 1's
arising in the 3-ary decomposition of $v(K)$.
\end{bem}

Now, for each tree $\phi$ with at least 2 leaves,
set $e(\phi):=(-1)^{w_\phi(1)+\ldots +w_\phi(n)}$.
Set $e(\tau):=1$.
For example, we have
$e(\beta)=-1$.
For the first tree in Example~\ref{extree}, we have $e(\phi)=-1$, and
for the second tree in Example~\ref{extree}, we have $e(\phi)=+1$.\\

Now, let $L$ be a graded module, $\phi$ an oriented tree with $n$
leaves and $B=(b_K)_{K\in V}$ a family of bilinear maps $L\ot L\nach L$.
Recursively, we want to define a multilinear map
$$\phi(B):L^{\ot n}\nach L.$$

\begin{enumerate}
\item
If $\phi$ has one leaf, i.e. $B$ is empty, we set $\phi(B):=\id$.
\item
If $\phi$ has only two leaves, i.e. $V=\{K_0\}$,
for a bilinear map $b_0: L\ot L\nach L$, we set
$\phi(b_0):=b_0$.
\item
If $\phi^{-1}(K_0)$ contains exactly one element, say $K_1$, and
$\pi_{K_0}(K_1)=1$, we set
$$\phi(B):=b_0\circ(\phi|_{K_1}((b_K)_{K\in V\setminus K_0})\ot 1).$$
\item
If $\phi^{-1}(K_0)$ contains exactly one element, say $K_1$, and
$\pi_{K_0}(K_1)=2$, we set
$$\phi(B):=b_0\circ(1\ot\phi|_{K_1}((b_K)_{K\in V\setminus K_0})).$$
\item
If $\phi^{-1}(K_0)=\{K_1,K_2\}$ with $\phi_{K_0}(K_1)=1$ and 
$\phi_{K_0}(K_2)=2$, we set
$$\phi(B):=b_0\circ(\phi|_{K_1}((b_K)_{K\in V_1})
\ot\phi|_{K_2}((b_K)_{K\in V_2})).$$
Here, $V_1$ denotes the ramification set of $\phi|_{K_1}$ and
$V_2$ the ramification set of $\phi|_{K_2}$.
\end{enumerate}

\subsection{Operations on trees}

\paragraph{Addition}
Let $(\phi,V,\pi)$ and $(\phi',V',\pi')$ be oriented trees with 
disjoint ramification sets. Let $R$ be a point in neither one 
of them.
Set $V'':=V\cup V'\cup \{R\}$. We define a map 
$\psi:V''\setminus R\nach V''$
by $\psi|_{V\setminus K_0}:=\phi$, $\psi|_{V'\setminus K'_0}:=\phi'$
and $\psi(K_0):=\psi(K'_0):=R$.\\

There is a family $(\pi''_K)_{K\in W}$ of inclusions
$\pi''_K:\psi^{-1}(K)\nach\{0,1\}$ with
$\pi''_K=\pi_K$, for $K\in V$, $\pi''_K=\pi'_K$ for $K\in V'$ and
$\pi''_R(K_0)=0$ and $\pi''_R(K'_0)=1$.
Now, we set
$$(\phi,V,\pi)+(\phi',V',\pi'):=(\psi,V'',\pi'').$$
For example, we have $\tau+\tau=\beta$. Each tree can be reconstructed
by addition out of copies of $\tau$.
It is obvious how to define the addition of non-oriented trees.
The addition of oriented trees is not commutative. The addition
of non-oriented trees is commutative.

\paragraph{Subtraction}
Let $(\phi,V)$ be a tree with $n$ leaves and $K\in V$.
Let $l$ be the number of leaves of $\phi|_K$. The definition
of a tree $\phi-\phi|_K$ with $n-l+1$ leaves is quite obvious.

\paragraph{Composition}
Let $(\phi,V,\pi)$ be an oriented tree with $n$ leaves and
let
$(\psi^{(1)},V^{(1)},\pi^{(1)}),$
$\ldots, (\psi^{(n)},V^{(n)},\pi^{(n)})$
be oriented trees. Let $W$ be the disjoint union of $V$ and all 
$V^{(i)}$.
For $K\in V$ set $n(K):=2-|\phi^{-1}(K)|$. (This is the number of leaves
belonging to $K$.) Let $K_1<\ldots<K_l$ all elements $K$ of $V$ with
$n(K)>0$. We define a map $\Phi:W\setminus K_0\nach W$ as follows:
For $K\in V\setminus K_0$, set $\Phi(K):=\phi(K)$.
For $K\in V^{(i)}\setminus K^{(i)}_0$, set $\Phi(K):=\psi^{(i)}(K)$,
and define the values of $\Phi$ on the $K^{(i)}$ by
$$(\Phi(K^{(1)}_0)\ddd \Phi(K^{(n)}_0)):=
(\underbrace{K_1\ddd K_1}_{n_1\text{ times}},\ldots,
\underbrace{K_l\ddd K_l}_{n_l\text{ times}}).$$
Then, $(\Phi,W)$ is a tree with a canonical orientation $\pi'$, 
given as follows:
For each $i$, $K\in V^{(i)}$ and $K'\in\Phi^{-1}(K)$, we set
$\pi'(K'):=\pi^{(i)}(K')$. For $K\in V$ and $K'\in\Phi^{-1}(K)\cap V$,
we set $\pi'(K'):=\pi(K')$. It remains to define $\pi'_{K_i}$ on elements
of $\Phi^{-1}(K_i)\setminus V$, for $i=1\ddd l$. So, if $n(K_i)$ equals
$2$, then $\Phi^{-1}(K_i)\setminus V$ has two elements, say
 $K^{(j)}_0$ and $K^{(k)}_0$ with $j<k$. Set $\Phi_{K_i}(K^{(j)}_0):=1$
and $\Phi_{K_i}(K^{(k)}_0):=2$.
If $n(K_i)$ equals
$1$, then $\Phi^{-1}(K_i)$ has one element in $V$, say $K$ and one 
element which is not in $V$, say $K'$.
Set $\Phi_{K_i}(K'):=1$ if $\phi_{K_i}(K)=2$ and
$\Phi_{K_i}(K'):=2$ if $\phi_{K_i}(K)=1$.\\

We will denote this decomposition by
$\Phi=\phi\circ(\psi^{(1)},\ldots ,\psi^{(n)})$.
The next lemma follows directly from the definitions:
\begin{lemma}\label{exponent}
In this situation, suppose that there is a family $B=(b_K)_{K\in W}$ of
bilinear maps $L\ot L\nach L$. Set $B^{(0)}:=(b_K)_{K\in V}$ and
$B^{(i)}:=(b_K)_{K\in V^{(i)}}$, for $i=1\ddd n$.
We have
$$\phi\circ(\psi^{(1)}\ddd\psi^{(n)})(B)=
(-1)^{\text{exponent} }
\phi(B^{(0)})\circ(\psi^{(1)}(B^{(1)})\ot\ldots\ot \psi^{(n)}(B^{(i)})),$$
where the exponent is the sum
$$(\sum_{K\in V^{(1)}}b_K)(\sum_{K\in V}^{V>1}b_K)+\ldots +
(\sum_{K\in V^{(n-1)}}b_K)(\sum_{K\in V}^{V>n-1}b_K).$$
We remind that $V>i$ means that the value $v(V)$ is greater
than the value $v(i)$ of the $i$-th leaf of $\phi$.
\end{lemma}

\section{Decomposition of differential graded Lie algebras}

\index{$L=(L,d,\lie)$}
Let $L=(L,d,\lie)$ be a differential graded Lie algebra, where the 
differential $d$ is of degree $+1$. Suppose that there is a splitting 
$\eta$, i.e. a map of degree $-1$ such that $d\eta d=d$. Furthermore, 
suppose that $\eta^2=0$ and $\eta d\eta=\eta$. When we use a Lie 
bracket on $\Hom(L,L)$, we mean the graded commutator. 
Observe that $[d,\eta]$ is a projector, and we can identify
the cohomology $H(L)$ with the graded submodule
$H:=\Fix(1-[d,\eta])$ of $L$. Set
$F$ to be the $d$-invariant submodule $\Fix([d,\eta])$ of $L$. 
We can consider the pair $(F,d)$ as an $L_\infty$-algebra.
As DG module, we have a decomposition
$(L,d)\isom (H,0)\oplus (F,d)$.
The aim of this section is to 
construct a minimal $L_\infty$-structure
on $H$ such that the decomposition $L\isom H\oplus F$
still holds in the category of $L_\infty$-algebras.
We begin with the construction of an $L_\infty$-algebra
map $g:F\to L$.\\

In the next theorem, for a tree $\phi\in\Ot(n)$,
we write $\Phi$ for the multilinear map
$\phi(\lie\ddd\lie):L^{\ot n}\nach L$.

\begin{satz}\label{ftol}
The following graded anti-symmetric maps
$g_n:F^{\ot n}\nach L$ of degree $1-n$ define a morphism
$g:F\nach L$ of $L_\infty$-algebras:
\begin{align*}
g_1:=&\text{inclusion}\\
g_2:=&\frac{-1}{2}\cdot[\eta\cdot,\cdot]\circ\alpha_2\\
    &\vdots\\
g_n:=&(\frac{-1}{2})^{n-1}\sum_{\phi\in\Ot(n)}\Phi\circ(\eta^{n-1}\ot 1)\circ
\alpha_n
\end{align*}
\end{satz}
\begin{proof}

\begin{align*}
&(-1)^{n-1}g_n\circ(d\ot 1^{\ot n-1})\circ\alpha_{1,n}=&\\[3pt]
&(\frac{1}{2})^{n-1}\sum_{\phi\in\Ot(n)}\sum_{i+j=n-1}\Phi\circ
(\eta^{\ot n-1}\ot 1)\circ(1^{\ot i}\ot d\ot 1^{\ot j})
\circ\alpha_n=&\\[3pt]
&(\frac{-1}{2})^{n-1}\sum_{\phi\in\Ot(n)}\sum_{i+j=n-2}(-1)^{i-1}
\cdot\Phi\circ
(\eta^{\ot i}\ot 1\ot\eta^{\ot j}\ot 1)\circ\alpha_n+&\\[3pt]
&(\frac{-1}{2})^{n-1}\sum_{\phi\in\Ot(n)}\sum_{i+j=n-1}
\Phi\circ(1^{\ot i}\ot d\ot 1^{\ot j})\circ
(\eta^{\ot n-1}\ot 1)\circ\alpha_n=&\\[3pt]
&(\frac{-1}{2})^{n-1}\sum_{i+j=n}\!\!\!(-1)^i
\!\!\!\sum_{\phi_1\in\Ot(i)}
\sum_{\phi_2\in\Ot(j)}[\Phi_1\circ(\eta^{i-1}\ot 1)\circ\alpha_i,
\Phi_2\circ(\eta^{j-1}\ot 1)\circ\alpha_j]
\circ\alpha_{i,n}+&\\[3pt]
&(\frac{-1}{2})^{n-1}\sum_{\phi\in\Ot(n)}d\circ
\Phi\circ(\eta^{\ot n-1}\ot 1)\circ\alpha_n=&\\[3pt]
&\frac{-1}{2}\sum_{i+j=n}(-1)^i[g_i,g_j]\circ\alpha_{i,n}+d\; g_n.
\end{align*}
Thus, condition (\ref{beding}) is verified.
\end{proof}

Remark that in the proof of Theorem~\ref{ftol},
we didn't make use of the Jacobi identity.
However, the Jacobi identity effects that
many of the terms in the decomposition of
$g_n$ annihilate each other. Before using
the formulas in a computer algorithm, one has to
study this more closely.\\

Now we come to the slightly more difficult construction
of an $L_\infty$-algebra structure
$\mu_\ast$ on $H:=H(L,d)$ with $\mu_1=0$
and of an $L_\infty$-quasi-isomorphism $f:H\to L$.
In the $A_\infty$-context, the existence of 
an $A_\infty$-structure on the cohomology of a DG algebra $A$
had already be shown (for $A$ connected) by Kadeishvili
\cite{Kade}, Gugenheim/ Stasheff \cite{GuSt} and
(in the general case) by Merkulov \cite{MerkK}. Merkulov gives a
recursion formula for construction of the higher products.
A similar recursion formula for the $L_\infty$-case can be found in
the article \cite{Hust} of Huebschmann and Stasheff. 
Kontsevich and Soibelman \cite{KontS} rewrote the
higher terms obtained by Merkulov's construction ($A_\infty$-case)
in terms of decorated trees. Their formulas are still
recursion formulas.\\

We have to make some preparations. First of all, there is
the following simple but important lemma:

\begin{lemma}\label{sechszudrei}
Let $n\geq 3$ be a natural number. There is a 1:1-correspondence
between triples $(\Phi,K,\sigma)$, where $\Phi=(\Phi,V,\pi)$
is an oriented tree with $n$ leaves, K is a ramification in $V$,
$\sig$ a permutation in $\Sigma_n$ and 6-tuples
$(k,\phi,\psi,\rho,\gamma,\delta)$, where $k$ is a natural number
with $2\leq k\leq n-1$, $\phi$ is a tree in $\Ot(k)$,
$\psi$ is a tree in $\Ot(l)$ where $l:=n+1-k$, $\rho$ is a
shuffle in $\Sh(k,n)$ and $\gamma\in\Sigma_l$, $\delta\in\Sigma_k$
are permutations.\\

{\footnotesize
\setlength{\unitlength}{0.5cm}
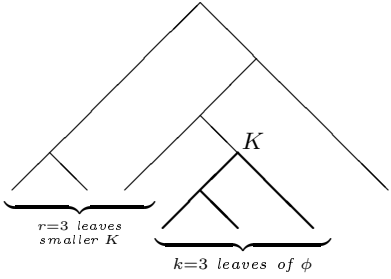
\hspace{-1.1cm}
\begin{picture}(12,9)
\put(7.1,4.1){$K$}
\put(0.8,2.8){$\underbrace{\quad\quad\quad\quad\quad\quad}_{r=3
\text{ leaves} \atop \text{smaller }K}$} 
\put(4.8,1.8){$\underbrace{\quad\quad\quad\quad\quad\quad\quad}_{k=3
\text{ leaves of }\phi}$} 
\put(1,3){\line(1,1){5}}
\put(2,4){\line(1,-1){1}}
\put(4,3){\line(1,1){3.5}}

\put(6,8){\line(1,-1){5}}
\put(6,5){\line(1,-1){1}}
\thicklines

\put(5,2){\line(1,1){2}}
\put(6,3){\line(1,-1){1}}
\put(7,4){\line(1,-1){2}}
\put(13,5){\parbox{8cm}{\textbf{Example:}\\ 
The fine lines represent the tree $\psi$
and the fat lines the tree $\phi$.\\ 
In the sequel, the first $r$ leaves of $\Phi$ will be
associated to the indexes $\sig(1)=\rho(\gamma(1)+k-1),\ldots
,\sig(r)=\rho(\gamma(r)+k-1)$, the following $k$ leaves
to the indexes $\sig(r+1)=\rho(\delta(1)),\ldots,
\sig(r+k)=\rho(\delta(k))$ and the remaining leaves 
to the indexes $\sig(r+k+1)=\rho(\gamma(r+2)+k-1),\ldots,
\sig(n)=\rho(\gamma(l)+k-1)$.}} 

\end{picture}}

To the triple $(\Phi,K,\sig)$, we associate the following
data:
Set $k$ to be the number of leaves of $\Phi|_K$, $\phi:=\Phi|_K$,
$\psi:=\Phi-\phi$. Let $r$ be the number of leaves $F$ of $\Phi$
with $F<K$. The shuffle $\rho$ is chosen in such a way that
$\{\rho(1),\ldots,\rho(k)\}=\{\sig(r+1),\ldots,\sig(r+k)\}$.
The permutation $\delta$ is defined by 
$\delta(i):=\rho^{-1}(\sig(r+i))$, for
$i=1,\ldots,k$ and $\gamma$ is defined in the following way:

$$\gamma(i):=\left\{ \begin{array}{r@{\quad\text{ for }\quad}l}
\rho^{-1}(\sig(i))-k+1 & i=1,\ldots,r\\
1 & i=r+1\\
\rho^{-1}(\sig(i+k-1))-k+1 & i=r+2,\ldots,l
\end{array}\right. $$ 

In the other direction, to the 6-tuple 
$(k,\phi,\psi,\rho,\gamma,\delta)$,
we associate the following data:
Set $r:=\gamma^{-1}(1)-1$. Then $\Phi$ is the composition
$$\Phi=\psi\circ(\underbrace{\tau,\ldots,\tau}_{r\text{ times}},\phi,
\tau,\ldots,\tau),$$
where $\tau$ again stands for the tree with one leaf.
The ramification $K$ is the root of $\phi$, considered as 
ramification of
$\Phi$ and $\sig$ is given by

$$\sig(i):=\left\{ \begin{array}{r@{\quad\text{ for }\quad}l}
\rho(\gamma(i)+k-1) & i=1,\ldots,r\\
\rho(\delta(i-r)) & i=r+1,\ldots,r+k\\
\rho(\gamma(i-(k-1))+k-1) & i=r+k+1,\ldots,n.
\end{array}\right. $$ 
\end{lemma}

Now suppose that such corresponding tuples $(\Phi,\hat{K},\sig)$
and $(k,\phi,\psi,\rho,\gamma,\delta)$ are given.
Let $V'$ be the ramification set of $\psi$ and $V''$ the ramification
set of $\phi$. Then $V:=V'\cup V''$ is the ramification set of
$\Phi$. Again, set $r:=\gamma^{-1}(1)-1$. 
Remark that the ordering on $V$ depends on $\gamma$.
We define a permutation $\tilde{\gamma}\in\Sigma_{l-1}$ by
  
$$\tilde{\gamma}(i):=\left\{ 
\begin{array}{r@{\quad\text{ for }\quad}l}
\gamma(i)-1 & i=1,\ldots, r\\
\gamma(i+1)-1 & i=r+1,\ldots,l-1.
\end{array}\right. $$ 

\begin{lemma}\label{wiwo}
We keep all notation from above. Let
$B=(b_K)_{K\in V}$ be a family of homogeneous bilinear forms
$L\ot L\nach L$. Denote the subfamilies $(b_K)_{K\in V'}$ and
$(b_K)_{K\in V''}$ by $B'$ and $B''$. Set $W$ to be the set
of all ramifications $K\in V$ such that $K>\hat{K}$. Then we have
\begin{multline*}
\psi(B')\circ\gamma\circ(\phi(B'')\circ\delta\ot 1\ot\ldots\ot 1)
\circ\rho\\[3mm]
=(-1)^{r+rk}\psi(B')\circ(\underbrace{1\ot\ldots\ot 1}_{r\text{ times}}
\ot\phi(B'')\ot 1\ot\ldots\ot 1)\circ\sig\\
=(-1)^{r+rk+\sum_{K\in W}b_k\cdot B''}\Phi(B)\circ\sig.
\end{multline*}
\end{lemma}

\begin{proof}
Let $a_1,\ldots,a_n$
be homogeneous elements of $L$. 
We get
\begin{multline*}
(\psi(B')\circ\gamma\circ(\phi(B'')\circ\delta\ot 1\ot\ldots\ot 1)
\circ\rho)(a_1\ot\ldots\ot a_n)=\\[3mm]
\chi(\rho,a_1,\ldots,a_n)\chi(\delta,a_{\rho(1)},
\ldots,a_{\rho(k)})\cdot\\[2mm]
(\psi(B')\circ\gamma)(
\underbrace{\phi(B'')(a_{\rho(\delta(1))},\ldots,a_{\rho(\delta(k))})}_{=:u_1}
\ot\underbrace{a_{\rho(k+1)}}_{u_2}\ot\ldots\ot
\underbrace{a_{\rho(n)}}_{u_l})=\\
\chi(\rho,a_1,\ldots,a_n)\chi(\delta,a_{\rho(1)},\ldots,a_{\rho(k)})
\chi(\gamma,u_1,\ldots,u_l)\psi(B')(u_{\gamma(1)},\ldots,u_{\gamma(l)}).
\end{multline*}
Using the following three formulas

\begin{align*}
\chi(\sig,a_1,\ldots,a_n)&=(-1)^{kr+(a_{\sig(1)}+\ldots+a_{\sig(r)})
(a_{\rho(1)}+\ldots+a_{\rho(k)})}\chi(\rho,a_1,\ldots,a_n)\cdot\\[2mm]
&\chi(\delta,a_{\rho(1)},\ldots,a_{\rho(k)})
\chi(\tilde{\gamma},a_{\rho(k+1)},\ldots,a_{\rho(n)}),\\[2mm]
u_{\gamma(1)}\ot\ldots\ot u_{\gamma(l)}&=
(-1)^{B''(a_{\sig(1)}+\ldots+a_{\sig(r)})}\cdot\\[2mm]
&(1\ot\ldots\ot 1\ot\phi(B'')\ot 1\ldots\ot 1)
(a_{\sig(1)}\ot\ldots\ot a_{\sig(n)}),\\[2mm]
\chi(\gamma,u_1,\ldots,u_l)&=
(-1)^{r+u_1(u_{\gamma(1)}+\ldots+ u_{\gamma(r)})}
\chi(\tilde{\gamma},u_2,\ldots,u_l),
\end{align*}
this expression is just
$$(-1)^{kr+r}\chi(\sig,a_1,\ldots,a_n)\psi(B')
((1\ot\ldots\ot 1\ot\phi(B'')\ot 1\ot\ldots\ot 1)
(a_{\sig(1)}\ot\ldots\ot a_{\sig(n)})).$$
The second equality of this Lemma is just a special case of
Lemma~\ref{exponent}.
\end{proof}

We turn to the construction of an $L_\infty$-structure
on $H(L)$.

\begin{prop}
The map
$[d,\eta]=d\eta+\eta d$ is a projection, i.e. $[d,\eta]^2=[d,\eta]$. 
And $H:=\Kern[d,\eta]$ is as module isomorphic to $H(L)$. 
Remark that under the assumption of the beginning of this section,
we have
$$H=\Kern d\cap\Kern\eta.$$
The bracket 
on $L$ induces a Lie-bracket on $H(L)$ and the induced bracket on $H$ 
(via the isomorphism $H\nach H(L)$) is just given by 
$(1-d\eta)\lie=(1-[d,\eta])\lie$.
\end{prop}

For simplicity, we set $g:=\eta\lie$.

\index{$\mu_\ast$}
\begin{satz}~\label{trick}
The following graded anti-symmetric maps $\mu_n:H^{\ot n}\nach H$ 
of degree $2-n$ define
the structure of an $L_\infty$-algebra on $H$:
\begin{align*}
\mu_1:&=0\\
\mu_2:&=(1-d\eta)\lie\\
&\vdots\\
\mu_n:&=(\frac{-1}{2})^{n-1}\sum_{\phi\in\Ot_n}e(\phi)\phi((1-[d,\eta])\lie,g
\ddd g)\circ\alpha_n
\end{align*} 
Here, the sum is taken over all trees $\phi$ with $n$ leaves
and $\phi((1-[d,\eta])\lie,g\ddd g)$ is the $n$-linear form obtained
by assigning the bilinear form $(1-[d,\eta])\lie$ to the root of the tree
$\phi$ and the bilinear form $g$ to each other ramification.
The sign $e(\phi)$ is defined in Section~\ref{trees}.
and $\alpha_n$ is the anti-symmetrization map.
\end{satz}
\begin{proof}
We must show that
\begin{equation}
\sum_{k+l=n+1}(-1)^{k(l-1)}\mu_l\circ(\mu_k\ot1\ot\ldots\ot 1)
\circ\alpha_{k,n}=0.
\end{equation}
Up to the factor $(-1)^{n-1}$ this sum has the form
\begin{multline}\label{nonc}
\sum_k
\sum_{\phi,\psi}\sum_{\rho,\gamma,\delta}(-1)^{k(l-1)}e(\phi)e(\psi)
\psi((1-[d,\eta])\lie,g,\ldots,g)\circ\gamma\circ\\
\circ
(\phi((1-[d,\eta])\lie,g,\ldots,g)\circ\delta\ot1\ot\ldots\ot1)\circ\rho,
\end{multline}
where $k$ ranges from from $2$ to $n-1$, $l=n+1-k$,
$\phi$ and $\psi$ vary in $\Ot(k)$ and $\Ot(l)$, $\rho$ in
$\Sh(k,n)$, $\gamma$ and $\delta$ in $\Sigma_l$ and $\Sigma_k$.
For corresponding tuples
$(k,\phi,\psi,\rho,\gamma,\delta)$ and $(\Phi,\hat{K},\sig)$
as in Lemma~\ref{sechszudrei}, we denote
as usual $r:=\gamma^{-1}(1)-1$ and by $t$ the number of ramifications
of $\psi$, greater than $r+1$.
Using
\begin{eqnarray*}
e(\Phi)&=&(-1)^{w_{\Phi}(\hat{K})(k-1)}e(\phi)e(\psi)\\
w_{\Phi}(\hat{K})&=&l -1-r-t
\end{eqnarray*} 
and Lemma~\ref{wiwo}, the expression (\ref{nonc}) can be expressed as
$$\sum_{\Phi\in\Ot(n)}\sum_{\hat{K}\in V\setminus K_0}e(\Phi)
(-1)^{r+w_\Phi(\hat{K})}\Phi(B)\circ\alpha_n,$$
where $B=(B_K)_{K\in V}$ is the family with 
$b_{K_0}=b_{\hat{K}}=(1-[d,\eta])\lie$ and $b_k=\eta\lie$
for $K\neq K_0,\hat{K}$.
To show that the last term is zero, it is enough to show
the following two conditions:
\begin{equation}
\sum_{\Phi\in\Ot(n)}\sum_{K\in V\setminus K_0}
\sum_{\sig\in\Sigma_n} (-1)^{r+w_\Phi(K)}e(\Phi)
\Phi((1-[d,\eta])\lie,g,\ldots,g,\underbrace{\lie}_{\text{pos. }K},
g,\ldots,g)\circ\sigma=0.
\end{equation}
For each tree $\Phi$, we have
\begin{equation}\label{nomal}
\sum_{K\in V\setminus K_0}
(-1)^{r+w_\Phi(K)} e(\Phi)
\Phi((1-[d,\eta])\lie,g,\ldots,g,
\underbrace{[d,\eta]\lie}_{\text{position }K},
g,\ldots,g)\circ\sigma=0.
\end{equation}

The first condition follows by the Jacobi-identity and an easy
combinatorial argument. In equation (\ref{nomal}) the term
annihilate each other since the differential d trickles down
the branches of $\Phi$:\\ 

\textbf{Initiation of the trickling:}
Suppose that $\Phi^{-1}(K_0)$ contains an element
$K'$ with $\pi_{K_0}(K')=1$. We have the following picture:
 
\setlength{\unitlength}{0.5cm}
{\footnotesize
\begin{picture}(23,5)
\put(0,2){$0$}
\put(1,2){$=$}
\put(3,2){\line(1,1){1}}
\put(4,3){\line(1,-1){1}}
\put(3,3.2){$(1-[d,\eta])d\lie$}
\put(2,1.2){$\eta\lie$}
\put(5,1.2){$\eta\lie$}

\put(7,2){$=$}

\put(10,2){\line(1,1){1}}
\put(11,3){\line(1,-1){1}}
\put(10,3.2){$(1-[d,\eta])\lie$}
\put(9,1.2){$d\eta\lie$}
\put(12,1.2){$\eta\lie$}

\put(14,2){$+(-1)^{\sharp\text{ramific. of }\Phi|_{K'}}$}
 
\put(20,2){\line(1,1){1}}
\put(21,3){\line(1,-1){1}}
\put(20,3.2){$(1-[d,\eta])\lie$}
\put(19,1.2){$\eta\lie$}
\put(22,1.2){$d\eta\lie$}
\end{picture}
}

\noindent
Here, we only have drawn the top of the tree $\Phi$
for the case where $\Phi^{-1}(K_0)$ consists of two
elements $K',K''$ and the corresponding bilinear
forms. It is quite obvious how this goes  when
$\Phi^{-1}(K_0)$ has only one element, since
$d|_H=0$.\\

\textbf{Going-on of the trickling} at a ramification
$K\in V$: 
\noindent
We illustrate the case, where $\Phi^{-1}(K)$ has two elements $K',K''$
with $\pi_K(K')=1$.

{\footnotesize
\begin{picture}(23,5)

\put(2,2){\line(1,1){1}}
\put(3,3){\line(1,-1){1}}
\put(2,3.2){$\eta d\lie$}
\put(1,1.2){$\eta\lie$}
\put(4,1.2){$\eta\lie$}

\put(6,2){$-$}

\put(9,2){\line(1,1){1}}
\put(10,3){\line(1,-1){1}}
\put(9.5,3.2){$\eta\lie$}
\put(8,1.2){$d\eta\lie$}
\put(11,1.2){$\eta\lie$}

\put(12,2){$-(-1)^{\sharp\text{ramific. of }\Phi|_{K'}}$}

\put(18,2){\line(1,1){1}}
\put(19,3){\line(1,-1){1}}
\put(18.5,3.2){$\eta\lie$}
\put(17,1.2){$\eta\lie$}
\put(20,1.2){$d\eta\lie$}
\put(21,2){$=\quad 0$}
\end{picture}
}

\noindent
Iterating the trickling down to the leaves and
using $d|_H=0$, we see that all terms in the sum are annihilated.
\end{proof}

\begin{bem}
The restriction of $\mu$ defines a formal map
$H^1\nach H^2$. One can show that this is just the
Kuranishi map as defined in \cite{Kosa}.
\end{bem}

\begin{satz}\label{trick2}
The following anti-symmetric maps $f_n:H^{\ot n}\nach L$ of degree
$1-n$ define an $L_\infty$-equivalence $H\nach L$ (i.e. an 
$L_\infty$-quasi-isomorphism).
\begin{align*}
f_1:&=\text{inclusion}\\
f_2:&=-g\\
&\vdots\\
f_n:&=-(\frac{-1}{2})^{n-1}\sum_{\phi\in\Ot(n)}e(\phi)
\phi(g\ddd g)\circ\alpha_n.
\end{align*}
\end{satz}

\begin{proof}
For $n\geq 0$, we have to prove the equation
\begin{equation*}
df_n-\sum_{i+j=n}\frac{(-1)^i}{2}[f_i,f_j]\alpha_{i,n}=
\sum_{k+l=n+1}(-1)^{k(l-1)}f_l\circ(\mu_k\ot 1\ot\ldots\ot 1)
\circ\alpha_{k,n}
\end{equation*} 
For $l=1$, the right hand-side is just $\mu_n$. Since
$$df_n=(\frac{-1}{2})^{n-1}\sum_{\phi\in\Ot(n)}
e(\phi)\phi(-d\eta\lie,g\ddd g)\circ\alpha_n,$$
it is sufficient to show the following three identities:

\begin{equation}\label{eins}
-\sum_{i+j=n}\frac{(-1)^{i}}{2}[f_i,f_j]\alpha_{i,n}=
(\frac{-1}{2})^{n-1}\sum_{\phi_\in\Ot(n)}e(\phi)
\phi(\lie,g\ddd g)\circ\alpha_n
\end{equation}
\begin{equation}\label{zwei}
f_l\circ(\phi(\lie,g\ddd g)\circ\alpha_k\ot 1\ot\ldots\ot 1)
\circ\alpha_{k,n}=0 \text{ for }  l>1,k+l=n+1.
\end{equation}

\begin{align}\label{drei}
&(\frac{-1}{2})^{n-1}\sum_{\phi_\in\Ot(n)}e(\phi)
\phi(\eta d\lie,g\ddd g)\circ\alpha_n=\nonumber\\
-&\sum_{k+l=n+1 \atop l\geq 2}
(-1)^{k(l-1)}\sum_{\phi\in\Ot(k)}(\frac{-1}{2})^{k-1}e(\phi)\quad\cdot\\
&f_l\circ(\phi(d\eta+\eta d)\lie,g\ddd g)\circ\alpha_k\ot1\ot
\ldots\ot 1)
\circ\alpha_{k,n}.\nonumber
\end{align}

\noindent
\textbf{Proof of equation (\ref{drei}):}
The right hand-side of equation (\ref{drei}) is
\begin{multline*}
(\frac{-1}{2})^{n-1}
\sum_{k+l=n+1}^{k,l\geq 2}\sum_{\phi,\psi}\sum_{\gamma,\delta,\rho}
(-1)^{k(l-1)}e(\phi)e(\psi)\psi(g,\ldots,g)\circ\gamma\circ\\
(\phi([d,\eta]\lie,g,\ldots,g)\circ\delta\ot 1\ot\ldots\ot 1)
\circ\rho.
\end{multline*}
As in the proof of Theorem~\ref{trick},
this expression takes the form
$$(\frac{-1}{2})^{n-1}
\sum_{\Phi\in\Ot(n)}\sum_{\hat{K}\in V}\sum_{\sig\in\Sigma_n}
(-1)^{r+w_\Phi(\hat{K})}e(\Phi)\Phi(B)\sig,$$
where $B=(b_K)_{K\in V}$ is
the family with $b_{\hat{K}}=[d,\eta]\lie$ and
$b_K=\eta\lie$ for $K\neq\hat{K}$.\\

\noindent
Hence to show equation (\ref{drei}), it is enough to show that 
for each tree $\phi$, we have
$$\Phi(\eta d\lie,g,\ldots,g)=\sum_{K\in V\setminus K_0}
(-1)^{r+w_\Phi(K)}\Phi(B).$$
This is true by the same trickling argument as in 
Theorem~\ref{trick}.\\

\noindent
\textbf{Proof of equation (\ref{zwei}):}
This is again the Jacobi-identity and some combinatorics.\\

\noindent
\textbf{Proof of equation (\ref{eins}):}  
\begin{multline*}
\sum_{i+j=n}\frac{(-1)^i}{2}\lie\circ(f_i\ot f_j)\circ\alpha_{i,n}=\\
=\sum_{i+j=n}(-\frac{1}{2})^{i-1+j-1}
\sum_{\phi\in\Ot(i),\psi\in\Ot(j)}\frac{(-1)^i}{2}
e(\phi)e(\psi)(\phi+\psi)(\lie,g,\ldots,g)\circ\alpha_n\\
=-(-\frac{1}{2})^{n-1}\sum_{\Phi\in\Ot(n)}
e(\Phi)\Phi(\lie,g,\ldots,g)\circ\alpha_n.
\end{multline*}
\end{proof}

Almost in the same manner, one can realize
the following: If $L=(L,d,\lie)$ is a DG Lie algebra
and 
$f_1:(M,d^M)\nach (L,d)$ a homotopy equivalence between DG modules,
then there is an $L_\infty$-algebra structure $\mu_\ast$ on
$M$ with $\mu_1=d^M$ and an $L_\infty$-quasi-isomorphism
$f:(M,\mu_\ast)\nach(L,d,\lie)$, extending $f_1$.
In other words: An up to homotopy differential graded
Lie algebra is an $L_\infty$-algebra.
In the $A_\infty$-context, this was already shown
by Markl \cite{Markl}.
Applying Proposition~\ref{isom} to the direct sum
$f\oplus g$, we get the 
decomposition theorem
for differential graded Lie algebras: 

\begin{satz}\label{decomposition}
The direct sum
$$f\oplus g:(H,\mu)\oplus (F,d)\nach (L,d,\lie)$$
is an isomorphism of $L_\infty$-algebras.
\end{satz}

One can show less constructively but more generally
(see \cite{KontL}) 
that each $L_\infty$-algebra with splitting
is isomorphic to the direct sum of a 
minimal and a linear contractible one.

\begin{kor}
If $(L,d,\lie)$ and $(L',d',\lie)$ are DG Lie algebras
such that $(L,d)$  and  $(L',d')$ are split, then,
for each $L_\infty$-quasi-isomorphism $f:L\nach L'$, there
exists an $L_\infty$-morphism $g:L'\nach L$ such that
$f_1$ and $g_1$ are inverse maps on the homology. 
\end{kor}

\section{Palamodov DG manifolds}

\subsection{Palamodov spaces}

We want to study a class of DG manifolds $(M,Q^M)$,
for which the vectorfield $Q^M$ is not only formal but
satisfies a convergence condition. Thus, the underlying 
graded vector space $M$
should be a normed or at least a metric vector space. 
It is natural, for instance, to define
``Banach DG manifolds'', where $M$ is a graded Banach space
and the formal sum $\sum_n\frac{1}{n!}\tilde Q^M$
converges in a neighbourhood of $0\in M$.
But in view of applications in analytic
deformation theory, this is not general enough
since, for example the tangent complex of a complex space
is not Banach.
We have to work more generally with
Palamodov's \textit{RO-} (or in Cyrillic \textit{PO-}) 
\textit{spaces}, whose
topology is defined by a countable family of semi-norms. 
We shall call these spaces \textit{Palamodov spaces} instead of
\textit{RO-spaces}. We will define \textit{Palamodov DG manifolds} 
and correspondingly \textit{Palamodov $L_\infty$-algebras}. 
This generalises
Bingener and Kosarew's notion \cite{BiKo}
of \textit{graded PO Lie algebras}. 

\begin{defi}
A \textbf{Palamodov space} is a vectorspace $E$ together with
a family $(\norm_\lambda)_{0<\lambda<1}$ of semi-norms such that
$$\norm_\lambda\leq\norm_{\lambda'}\quad\text{ for }\quad\lambda<\lambda'.$$
A morphism $f:E\nach F$ of Palamodov spaces is a morphism
of vectorspaces together with global constants $\epsilon\in (0,1)$
and $C>0$ such that, for each element $e\in E$ and each 
$\lambda\in(0,1)$, 
we have
$$||f(e)||_\lambda\leq C\cdot||e||_\lambda.$$
\end{defi}

\begin{beisp}\label{holomorph}
If $E$ is the space of 
holomorphic functions on an open polydisk $P$
of poly-radius $r\in\RR_{>0}^N$ we
can define seminorms 
$$|f|_\lambda:=\sup\{|f(z)|:\;z\in\lambda\cdot P\},$$
for $0<\lam<1$. 
By Cauchy's integral
formula, we have the estimation
$$(\lam'-\lam)^2\cdot|\d f/\d x_i|_\lam\leq |f|_{\lam'},$$
for $\lam<\lam'<1$.
If we set $$||\sum_\nu a_\nu x^\nu||_\lam:=
\sum_\nu|a_\nu|\cdot|x^\nu|_\lam
=\sum_\nu|a_\nu|\cdot\lam^{|\nu|}\cdot r^\nu,$$
and $C:=\frac{1}{e\cdot r_i}$,
we get the better estimation
$$(\lam'-\lam)\cdot ||\d f/\d x_i||_\lam\leq\frac{C}{\epsilon}
\cdot||f||_{\lam'},$$
for
$1-\epsilon<\lam<\lam'<1$.
For the following reason, we usually use the Palamodov
structure given by $||\cdot||_\lam$ on $E$:
Let $R^0=\Gamma(P,\Oh_P)$ be the $\CC$-algebra of
holomorphic functions on the open polydisk $P$ and
$R=R^0[x_i]_{i\in I}$ a 
finitely generated free DG algebra over $R^0$ with 
$g(x_i)<0$, for all $i\in I$.
The graded vector space $L=\Der(R,R)$ is then
a sum of copies of $E$ and carries
the structure of a Palamodov DG Lie algebra (see 
the definition below).
\end{beisp}

\paragraph{Multilinear maps and power series}

In this section,
$E$ and $F$ are graded vector spaces. 
By definition, a \textbf{polynomial 
$E\to F$ of (polynomial) degree $p$
and linear degree i} is a map of the form
$x\mapsto\tilde\phi(x)$, where
$\phi:E^{\odot p}\to F$ is a \textit{graded symmetric} 
$p$-linear map (of linear degree $i$) 
and $\tilde\phi(x)=\phi(x\ddd x)$.
Following Bingener and Kosarew \cite{BiKo},
for Palamodov spaces $E$ and $F$ and a given $\epsilon\in (0,1)$, we
define the following pseudo-norms on the graded
vectorspace
of polynomial maps $u:E\to F$ of polynomial degree $p$:
$$||u||^{0,\epsilon}:=\sup\{||u(x)||_\lambda:
\;||x||_{\lambda}\leq 1\text{ and } 1-\epsilon\leq\lambda<1\}$$
$$||u||^{1,\epsilon}:=\sup\{(\lambda'-\lambda)^{p-1}||u(x)||_\lambda:
\;||x||_{\lambda'}\leq 1\text{ and } 1-\epsilon\leq\lambda<\lambda'<1\}$$
Furthermore, we define the following pseudo-norms
on the graded vectorspace of $p$-linear maps $\phi:E^{\ot p}\to F$.
$$|\phi|^{0,\epsilon}:=\sup\{
||\phi(x_1\ddd x_p)||_\lambda:
\;||x_i||_{\lambda}\leq 1\quad\forall i \text{ and } 
1-\epsilon\leq\lambda<1\}$$
$$|\phi|^{1,\epsilon}:=\sup\{(\lambda'-\lambda)^{p-1}
||\phi(x_1\ddd x_p)||_\lambda:
\;||x_i||_{\lambda'}\leq 1\quad\forall i \text{ and }
1-\epsilon\leq\lambda<\lambda'<1\}$$

\begin{bem}\label{untilde}
Suppose that $E$ and $F$ are graded Palamodov spaces. 
For each graded symmetric $p$-linear map
$\phi:E^{\odot p}\to F$, we
have
$$||\phi||^{1,\epsilon}\leq |\phi|^{1,\epsilon}\leq\frac{p^p}{p!}
||\tilde\phi||^{1,\epsilon}.$$
\end{bem}

For graded vector spaces $E$ and $F$, set
$F[[E]]^i$ to be the set of formal sums
$\sum_{p=0}^\infty u_p$, where $u_p:E\to F$
is a polynomial of degree $p$ and linear degree $i$.
By omission, we mean $i=0$ or another fixed $i\in\ZZ$. 
A formal map $u=\sum_{p=0}^\infty u_p$ in $F[[E]]$
is called \textbf{strictly convergent} (with respect to $\epsilon$), 
if there exists a $t>0$ such that
$$\sum_{p=0}^\infty||u_p||^{0,\epsilon}\cdot t^p<\infty.$$ 
A power series $u=\sum_{p=0}^\infty u_p$ in $F[[E]]$
is called \textbf{1-convergent} (with respect to $\epsilon$), if
there exists a $t>0$ such that
$$\sum_{p=0}^\infty||u_p||^{1,\epsilon}\cdot t^p<\infty.$$ 

Each strictly convergent formal map is 1-convergent.

\paragraph{Palamodov DG manifolds}

\begin{defi}\label{CODER}
For $i\in\{-max\ddd +max\}$, let $M^i$ be Palamodov spaces 
and $M:=\oplus M^i$. A coderivation
$Q\in\Coder^{\geq 0}(SM)$ on $SM$ is called
\textbf{1-convergent} (resp. \textbf{strictly convergent})
with respect to $\eps$, if the formal sum
$\tilde{Q}:=\sum_n\frac{1}{n!}\tilde{Q}_n$ in $M[[M]]^1$
is 1-convergent
(resp. strictly convergent) with respect to $\epsilon$.
A coderivation $Q$ is called \textbf{convergent}
with respect to $\eps$, if
for each component $Q_n\in\Mult_n(M,M)$,
we have $|Q_n|^{1,\epsilon}<\infty$ and if the formal
sum $\tilde{Q}:=\sum_n\frac{1}{n!}\tilde{Q}_n$, 
considered as formal map 
$M^{\geq 0}\nach M^{\geq 0}$ is strictly convergent 
with respect to $\epsilon$. Denote the set of all convergent
coderivations of $SM$ of degree $i$ by $\CODER^i(SM)$ and set
$\CODER(SM):=\oplus_{i\geq 0}\CODER^i(SM)$.
\end{defi}

We want to emphasise
that by graded symmetry,
the odd components of $m$ only give a linear contribution to
the term $\tilde{Q}(m)$.

\begin{prop}
The sets of 1-convergent, strictly convergent and
convergent coderivations, respectively  of $SM$
are graded Lie subalgebras of $\Coder^{\geq 0}(SM)$.
\end{prop}
\begin{proof}
We only prove the statement for 1-convergent coderivations.
We must show that, if $Q$ and $q$ are 1-convergent coderivations,
then $[Q,q]$ is 1-convergent. Fix $\lambda_0,\lambda_1$ such that
$1-\epsilon\leq\lambda_0<\lambda_1<1$. Set 
$\lambda:=\frac12(\lambda_1-\lambda_0)$.

\begin{equation*}
\begin{array}{c}
(\lam_1-\lam_0)^{n-1} || [Q,q]^\sim_n ||_{\lam_0} \leq\\[3pt]
2^{n-1}(\lam-\lam_0)^{n-1}\sum_{k+l=n+1}(
||Q_l(\tilde q_k(x),x\ddd x)||_{\lam_0}+
||q_l(\tilde Q_k(x),x\ddd x)||_{\lam_0})\leq\\[4pt]
2^{n-1}(\lam_1-\lam)^{k}\sum_{k+l=n+1}(
|Q_l|^{1,\eps}||\tilde q_k(x)||_\lam ||x||_\lam^{l-1}+
|q_l|^{1,\eps}||\tilde Q_k(x)||_\lam ||x||_\lam^{l-1})\leq\\[4pt]
2^{n-1}\sum_{k+l=n+1}(
|Q_l|^{1,\eps} ||q_k||^{1,\eps}+|q_l|^{1,\eps} ||Q_k||^{1,\eps})
||x||_{\lam_1}^n\leq\\[4pt]
\frac{2^{n-1}(n-1)^{n-1}}{(n-1)!}\sum
(||Q_l||^{1,\eps} ||q_k||^{1,\eps}+||q_l||^{1,\eps} ||Q_k||^{1,\eps})
||x||_{\lam_1}^n
\end{array}
\end{equation*}  
Thus,
\begin{equation*}
\begin{array}{c}
|| [Q,q]^\sim_n ||^{1,\eps}\leq
\frac{2^{n-1}(n-1)^{n-1}}{(n-1)!}\sum
(||Q_l||^{1,\eps} ||q_k||^{1,\eps}+||q_l||^{1,\eps} ||Q_k||^{1,\eps}).
\end{array}
\end{equation*}
By the Stirling formula, we have 
$\frac{2^{n-1}(n-1)^{n-1}}{(n-1)!}\leq (2e)^{n-1}$, hence
\begin{equation*}
\sum_n\frac{1}{n!}|| [Q,q]^\sim_n ||^{1,\eps}t^n\leq
2(\sum_k\frac{(2e)^k}{k!}||\tilde Q_k||^{1,\eps}t^k)
(\sum_l\frac{(2e)^l}{l!}||\tilde q_l||^{1,\eps}t^l),
\end{equation*}
which is finite, for small $t$.

\end{proof}

\begin{defi}
A formal DG manifold $(M,Q^M)$, where $M=\oplus_{i=-d}^d M^i$ is
a graded Palamodov space and $Q^M\in\CODER^1(SM)$ will be called
\textbf{Palamodov DG manifold}. If $M$ is Banach and the 
Palamodov structure
is 
given by $||\cdot||_\lam:=||\cdot||$, for each $\lam\in(0,1)$, 
then $(M,Q^M)$ will be called
\textbf{Banach DG manifold}.
A \textbf{morphism} $f:M\nach N$ of Palamodov DG manifolds is
a morphism $f$ of formal DG manifolds such that $||f_n||^{1,\eps}$
is finite, for each $n$, and such that the formal map
$\tilde f=\sum_n\frac{1}{n!}\tilde f_n$ is strictly 
convergent as a map
$M^{\geq 0}\nach M^{\geq 0}$.  
\end{defi}

The following
inverse function theorem shows that Palamodov
DG manifolds are well behaved:

\begin{prop}
If $f:M\to N$ is a morphism of Palamodov super manifolds
(resp. Palamodov DG manifolds) such that $f_1$
is an isomorphism of Palamodov spaces, then $f$
is an isomorphism.
\end{prop}
\begin{proof}
The proof of Proposition II.8.4 of \cite{BiKo}
applies directly to this situation. For convenience,
we indicate it, however. It uses the following
statement: Let $\sum_p\gamma_p t^p$ the Taylor
series at zero of the function $\frac 12-\frac 14\sqrt{1-t}$.
The $\gamma_p$ are all positive and have the property
that for any $p,q>0$ the inequality
$\sum_{i_1\kkk i_p=q}\gamma_{i_1}\ppp\gamma_{i_p}\leq\gamma_q$
holds. A power series $\sum_pa_pt^p$ with nonnegative
coefficients $a_p$ converges near zero, if and only if
there exist constants $C,R>0$ such that
$a_p\leq\gamma_p\cdot C\cdot R^p$, for all $p\geq 0$.\\

For $n\geq 1$, let $g_n:N^{\odot n}\to M$ be the map
defined as in the proof of Proposition~\ref{isom}.
The only thing to verify is that the formal
map $\sum_q\frac{1}{q!}\tilde g_q:N^{\geq 0}\to M^{\geq 0}$ 
strictly converges. Without restriction, we may assume
that $M$ and $N$ are positively graded.
By assumption, there are real numbers $C,R>0$ such that
$\frac{1}{p!}|f_p|^{0,\eps}\leq\gamma_p\cdot C\cdot R^p$, for each $p$.
Chose $C',R'>0$ such that
$$||g_1||^{0,\eps}\leq\gamma_1\cdot C'\cdot R'$$
and
$$(\sum_{p=2}^\infty\gamma_p(C'\cdot R)^{p-1})\cdot
||g_1||^{0,\eps}\cdot C\cdot R\leq 1.$$
We show by induction on $q$ that
$\frac{1}{q!}||\tilde g_q||^{0,\eps}\leq\gamma_q\cdot C'\cdot (R')^q$,
for each $q$. The case $q=1$ is already done. For $q\geq 2$, we
get:
$$\begin{array}{l}
\frac{1}{q!}||\tilde g_q||^{1,\eps}\leq\\[8pt]
\sum_{p=2}^n\sum_{i_1\kkk i_p=n}\frac{1}{p!I!}||g_1\circ f_p\circ
(g_{i_1}\ots g_{i_p})||^{0,\eps}\leq\\[8pt]
\sum_{p=2}^n\sum_{i_1\kkk i_p=n}\frac{1}{p!I!}||g_1||^{0,\eps}
|f_p|^{0,\eps}||g_{i_1}||^{0,\eps}\ppp||g_{i_p}||^{0,\eps}\leq\\[8pt]
\sum_{p=2}^n\sum_{i_1\kkk i_p=n}||g_1||^{0,\eps}\gamma_p C R^p
\gamma_{i_1}\ppp\gamma_{i_p} (C')^p(R')^q\leq\\[8pt]
(\sum_{p=2}^n\gamma_p(C' R)^{p-1})||g_1||^{0,\eps}C R C'
\gamma_q (R')^q\leq\\[8pt]
\gamma_q C' (R')^q.
\end{array}$$
\end{proof}

\subsection{Decomposition of Palamodov DG Lie algebras}

In this section, we will prove our main result, an
analytic version of the Decomposition Theorem~\ref{decomposition}.
To this aim,
suppose that $L=(L,d,\lie)$ is a Palamodov Lie algebra,
equipped with a splitting $\eta$ of $d$ such that each component
$\eta^i:L^i\nach L^{i-1}$ is a morphism of Palamodov spaces.
In general, it is very difficult to construct such a splitting.
For the case, where $L$ is the tangent complex of a complex
space, this was realized by Palamodov.
Generalisations of Palamodov's splitting construction 
can be found in \cite{BiKo}.
In the Banach situation, the construction is much easier.\\
 
Again, without loss of generality, we may assume
that $\eta^2=0$, $\eta d\eta=\eta$ and $d\eta d=d$.
From now on, $L$ stands for the quadruple $(L,d,\lie,\eta)$.
We define the $L_\infty$-structure on $H=\Fix(1-[d,\eta])$
an the $L_\infty$-morphism $f:H\nach L$
exactly in the same way as in Theorem~\ref{trick} and \ref{trick2}. 
Both, $H$ and the complement $F$ of $H$ in $L$ inherit the
graded Palamodov space structure of $L$.
First, we show that $H=(H,\mu_\ast)$ is a Palamodov $L_\infty$-algebra.

\begin{lemma}\label{schaetzen}
Let $E$ be a Palamodov space, $\phi\in\Ot(n)$ a binary tree with
$n$ leaves and $g_1\ddd g_{n-1}\in\Mult_2(E,E)$ bilinear forms.
For $\eps\in(0,1)$ and each $i=1\ddd n-1$, 
suppose that $b_i:=||g_i||^{0,\epsilon}$
and  $c_i:=||g_i||^{1,\epsilon}$ are finite. Then, we have 
the following inequalities:
$$\begin{array}{c}
||\phi(g_1\ddd g_{n-1})||^{0,\eps}\leq b_1\ppp b_{n-1}\\
||\phi(g_1\ddd g_{n-1})^\sim||^{1,\eps}\leq (n-1)^{n-1}c_1\ppp c_{n-1}
\end{array}$$
\end{lemma}
\begin{proof}
The first inequality is easy. 
We prove the second one by induction on $n$.
There is only one tree $\beta$ with two leaves and $\beta(g_1)=g_1$.
This proves the case $n=2$. Suppose that the statement is true, for
each $m\leq n-1$. We have to show that for fixed $\lam_0<\lam_1$
in $[1-\eps,1)$ and for each element $x\in E$, we have
\begin{equation}\label{toshow}
(\lam_1-\lam_0)^{n-1}||\phi(g_1\ddd g_{n-1})^\sim(x)||_{\lam_0}\leq
(n-1)^{n-1}c_1\ppp c_{n-1}\cdot ||x||_{\lam_1}^n.
\end{equation} 
We have $\phi(g_1\ddd g_{n-1})=g_1\circ\;\psi_1(g_2\ddd g_p)\ot
\psi_2(g_{p+1}\ddd g_{n-1})$, for some $p\leq n-2$ and trees
$\psi_1\in\Ot(p)$ and $\psi_2\in\Ot(q)$.
We will just write $\phi$ instead of $\phi(g_1\ddd g_{n-1})^\sim$,
$\psi_1$ instead of $\phi(g_2\ddd g_{k})^\sim$ and
$\psi_2$ instead of $\phi(g_{k+1}\ddd g_{n-1})^\sim$.
Set $\lam:=\lam_0+\frac{1}{n-1}(\lam_1-\lam_0)$. Then, we get
$$\begin{array}{c}
(\lam_1-\lam_0)^{n-1}\cdot||\phi(x)||_{\lam_0}=\\[3pt]
(n-1)^{n-1}(\lam-\lam_0)||g_1(\psi_1(x),
\psi_2(x))||_{\lam_0}\leq\\[3pt]
c_1(n-1)^{n-1}(\lam-\lam_0)^{n-2}||\psi_1(x)||_{\lam}
||\psi_2(x)||_{\lam}=\\[3pt]
c_1\frac{(n-1)^{n-1}}{(n-2)^{n-2}}(\lam_1-\lam)^{p-1}
\cdot||\psi_1(x)||_{\lam}\cdot(\lam_1-\lam)^{q-1}\cdot
||\psi_2(x)||_{\lam}\leq\\[3pt]
c_1\frac{(n-1)^{n-1}}{(n-2)^{n-2}}||\psi_1||^{1,\eps}
\cdot ||x||^p_{\lam_1}
\cdot ||\psi_2||^{1,\eps}\cdot ||x||^q_{\lam_1}.
\end{array}$$
By the induction hypothesis, the last term is lower or equal 
$$\begin{array}{c}
c_1\frac{(n-1)^{n-1}}{(n-2)^{n-2}}(p-1)^{p-1}c_2\ppp c_p
(q-1)^{q-1}c_{p+1}\ppp c_{n-1}\cdot||x||^n_{\lam_1}\leq\\[3pt]
(n-1)^{n-1}c_1\ppp c_{n-1}\cdot||x||^n_{\lam_1}.
\end{array}$$
\end{proof}

\begin{satz}
The $L_\infty$-algebra $(H,\mu_\ast)$ is Palamodov and the
$L_\infty$-maps $f:H\nach L$ 
and $g:F\nach L$ are morphism of Palamodov $L_\infty$-algebras.
\end{satz}
\begin{proof}
The proofs of all three statements go almost in the
same way, thus we only prove the first statement.
Set $c:=||\eta||^{0,\epsilon}$ and $k:=||\lie||^{1,\epsilon}$.
By Lemma~\ref{schaetzen},
we have
$$||\phi(\lie,\eta\lie\ddd\eta\lie)||^{1,\eps}\leq
(n-1)^{n-1}k^{n-1}c^{n-2},$$
hence 
\begin{align*}
||\tilde\mu_n||^{1,\epsilon}\leq\quad & n!(\frac{(n-1)}{2})^{n-1}
\#\Ot(n)\cdot k^{n-1}\cdot c^{n-2}\\
\leq\quad & n!(n-1)^{n-1}
\cdot (2k)^{n-1}\cdot c^{n-2}\quad<\quad\infty.
\end{align*}
The last step is by Lemma~\ref{otnum}. This proves the first
condition for the convergence of $\mu_\ast$. To prove
the second condition, set
$$\kappa:=|\text{(co)restriction of $\lie$ to $L^{\geq 1}\ot L^{\geq 1}$}
|^{0,\eps},$$ 
which is finite by the assumption that $L$ is Palamodov.
By Lemma~\ref{schaetzen} and \ref{otnum}, we get
\begin{align*}
||\tilde{\mu}_n||^{0,\eps}\quad\leq\quad 
\frac{n!}{2^{n-1}}\#\Ot(n)\kappa^{n-1}c^{n-1}\quad
=\quad 2^{n-1}n!\kappa^{n-1}c^{n-2}.
\end{align*}
Thus, for small $t>0$, the series
$\sum_n\frac{1}{n!}||\tilde{\mu}_n||^{0,\eps}\cdot t^n$ converges.
\end{proof}

Applying the inverse function theorem, we
finally get the main result:

\begin{kor}
If $(L,d,\lie,\eta)$ is a Palamodov DG Lie algebra
with splitting $\eta$, then 
we have a decomposition
$$(L,d,\lie)\isom(H,\mu_\ast)\oplus(F,d)$$ 
in the category of Palamodov $L_\infty$-algebras.
\end{kor}

\end{document}